\documentclass[10pt,a4paper]{amsart}
\usepackage{hyperref}
\hypersetup{hypertexnames = false, bookmarksdepth = 2, bookmarksopen = true, colorlinks, linkcolor = black, citecolor = black, urlcolor = black, pdfstartview={XYZ null null 1}}

\usepackage{amsmath,amssymb,amsthm}
\usepackage{booktabs}
\usepackage{enumitem}
\usepackage{fixltx2e}
\usepackage{mathtools}
\usepackage{tikz}
\usepackage[textsize=small]{todonotes}
\usepackage{xparse}
\usepackage{xspace}
\usepackage{mathtools}

\usepackage{mathrsfs}
\usepackage{microtype}
\numberwithin{equation}{section}
\def\mod{\operatorname{mod}}

\def\CH{\mathop{\mathrm{CH}}}

\def\Ext{\operatorname {Ext}}

\def\End{\operatorname {End}}

\def\End{\operatorname {End}}

\def\rk{\operatorname {rk}}

\def\Pic{\operatorname {Pic}}

\theoremstyle{definition}

\newtheorem{lemma}{Lemma}[section]
\newtheorem{proposition}[lemma]{Proposition}
\newtheorem{theorem}[lemma]{Theorem}
\newtheorem{corollary}[lemma]{Corollary}

\newtheorem{example}[lemma]{Example}
\newtheorem{definition}[lemma]{Definition}

\newtheorem{remark}[lemma]{Remark}

\newcommand\Bl{\ensuremath{\mathrm{Bl}}}
\newcommand\bounded{\ensuremath{\mathrm{b}}}
\newcommand\catmod{\ensuremath{\mathrm{mod}}} 
\newcommand\coh{\ensuremath{\mathrm{coh}}}

\newcommand\num{\ensuremath{\mathrm{num}}}

\newcommand\An{\ensuremath{\mathsf{A}}}
\newcommand\Dn{\ensuremath{\mathsf{D}}}
\newcommand\En{\ensuremath{\mathsf{E}}}
\newcommand\kronecker{\ensuremath{\mathsf{K}}}
\newcommand\Sn{\ensuremath{\mathsf{S}}}

\DeclareMathOperator\derived{\mathbf{D}}

\DeclareMathOperator\HHHH{HH}
\DeclareMathOperator\Hom{Hom}

\DeclareMathOperator\Kzero{K_0}
\DeclareMathOperator\Kzeronum{K_0^{num}}

\DeclareMathOperator\PProj{\mathbf{Proj}}
\DeclareMathOperator\RRRHom{\mathbf{R}Hom}
\DeclareMathOperator\serre{\mathbb{S}}

\mathchardef\mhyphen="2D
\newcommand\dash{\nobreakdash-\hspace{0pt}}

\newcounter{todocounter}
\DeclareDocumentCommand\addreference{g}{\stepcounter{todocounter}\todo[color = blue!30, fancyline]{\thetodocounter. Add reference\IfNoValueF{#1}{: #1}}\xspace}
\DeclareDocumentCommand\checkthis{g}{\stepcounter{todocounter}\todo[color = red!50, fancyline]{\thetodocounter. Check this\IfNoValueF{#1}{: #1}}\xspace}
\DeclareDocumentCommand\fixthis{g}{\stepcounter{todocounter}\todo[color = orange!50, fancyline]{\thetodocounter. Fix this\IfNoValueF{#1}{: #1}}\xspace}
\DeclareDocumentCommand\expand{g}{\stepcounter{todocounter}\todo[color = green!50, fancyline]{\thetodocounter. Expand\IfNoValueF{#1}{: #1}}\xspace}

\title{Embeddings of algebras in derived categories of surfaces}
\author{Pieter Belmans \and Theo Raedschelders}

\begin{document}

\begin{abstract}
  By a result of Orlov there always exists an embedding of the derived category of a finite-dimensional algebra of finite global dimension into the derived category of a high-dimensional smooth projective variety. In this article we give some restrictions on those algebras whose derived categories can be embedded into the bounded derived category of a smooth projective surface. This is then applied to obtain explicit results for hereditary algebras.
\end{abstract}

\maketitle

\newcommand\mscfootnote[1]{%
  \begingroup
  \renewcommand\thefootnote{}\footnote{#1}%
  \addtocounter{footnote}{-1}%
  \endgroup
}

\mscfootnote{\emph{2010 Mathematis Subject Classification.} Primary 14F05, 16E35; Secondary 18E30.}

\maketitle

\section{Introduction}
\label{section:introduction}
For an algebraically closed field $k$ of characteristic~$0$, consider a triangulated category of the form~$\derived^\bounded(\mod/A)$, for some finite-dimensional $k$\dash algebra~$A$ of finite global dimension. In two recent preprints, Orlov~\cite{orlov-nc-schemes,orlov-geometric-realization} showed that there always exists an admissible embedding
\begin{equation}
  \label{equation:embedding}
  \derived^\bounded(\catmod/A) \hookrightarrow \derived^\bounded(\coh/X),
\end{equation}
for some smooth projective variety~$X$.

This construction typically embeds~$A$ into a high-dimensional variety, and we consider the existence of an embedding in the derived category of a surface, as the case of curves is completely understood by the indecomposability result of Okawa showing there are no non-trivial embeddings \cite{okawa-sod-curve}. In particular we give two types of obstructions to the existence of an embedding~\eqref{equation:embedding} where~$X=S$ is a smooth projective surface. These take the form of conditions on the Euler form of the algebra~$A$. We summarise these as follows.

\begin{theorem}(see corollaries \ref{corollary:rank}, \ref{corollary:signature} below)
  \label{theorem:main-theorem}
  If an embedding
  \begin{equation}
    \derived^\bounded(\catmod/A) \hookrightarrow \derived^\bounded(\coh/S),
  \end{equation}
  exists, then~$\rk(\chi^-_{A}) \leq 2$ and~$\chi_{A}^+$ does not admit a~$3$\dash dimensional negative definite subspace.
\end{theorem}

It is not hard to see that the use of non-commutative motives and additive invariants doesn't yield any strong results. Hence to obtain these results one has to incorporate extra data coming from the Euler form on the triangulated categories, together with an understanding of the structure of the numerical Grothendieck groups of surfaces.

This type of result explains one aspect of the structure of semi-orthogonal components of the derived category of a surface. On one hand it is known that there exist (quasi-)phantoms for surfaces, which are components that are very hard to understand \cite{gorchinskiy-orlov-geometric-phantom-categories}. This article on the other hand studies certain ``easy'' components coming from derived categories of finite-dimensional algebras, in particular algebras arising from strong exceptional collections. One interesting result from these obstructions is the non-existence of an exceptional collection of~4~objects whose endomorphism algebra is isomorphic to~$k\An_4$.

Section~\ref{section:preliminaries} contains some preliminaries on triangulated categories and the structure of derived categories of finite-dimensional algebras and smooth projective varieties. In section~\ref{section:obstructions} we discuss the constraints that are imposed by an embedding \eqref{equation:embedding} by using noncommutative motives and incorporating the extra data. In section~\ref{section:applications-constructions} we discuss examples of finite-dimensional algebras that violate the conditions, and we give explicit embeddings for families of algebras that do satisfy the constraints. We also list some open questions on the structure of strong exceptional collections inside derived categories of surfaces raised by the obstructions and examples.

\paragraph{\textbf{Acknowledgements}}
We would like to thank Louis de Thanhoffer de V\"olcsey and Michel Van den Bergh for interesting discussions, David Ploog and Alexey Elagin for their comments on the first version, and Nathan Prabhu--Naik and Markus Perling for answering questions regarding toric geometry.

Both authors were supported by a Ph.D.\ fellowship of the Research Foundation---Flanders (FWO).

\section{Preliminaries}
\label{section:preliminaries}
Let~$\mathcal{T}$ denote a~$k$\dash linear triangulated category, where~$k$ denotes an algebraically closed field of characteristic~$0$. Based on~\cite{orlov-nc-schemes} we recall a couple of definitions and lemmas that we will require later on. Let~$\mathcal{N}$ denote a full triangulated subcategory of~$\mathcal{T}$.

\begin{definition}
  A full embedding~$i: \mathcal{N} \hookrightarrow \mathcal{T}$ is \emph{left} (respectively \emph{right}) \emph{admissible} if there is a left (respectively right) adjoint functor~$q\colon\mathcal{T} \to \mathcal{N}$ to~$i$. It is \emph{admissible} if~$i$ is both left and right admissible.
\end{definition}

Recall that~$\mathcal{N}^{\perp}$ denotes the right orthogonal to~$\mathcal{N}$: it is the full subcategory of~$\mathcal{T}$ consisting of objects~$M$ such that~$\Hom_{\mathcal{T}}(\mathcal{N},M)=0$. The left orthogonal is defined similarly and is denoted~${}^{\perp}\mathcal{N}$.

\begin{definition}
  The triangulated category~$\mathcal{T}$ has a \emph{semi-orthogonal decomposition}
  \begin{equation}
    \mathcal{T}=\langle \mathcal{N}_1,\ldots,\mathcal{N}_n \rangle
  \end{equation}
  for full triangulated subcategories~$\mathcal{N}_i$, if~$\mathcal{T}$ has an increasing filtration
  \begin{equation}
    0=\mathcal{T}_0 \subset \mathcal{T}_1 \subset \cdots \subset \mathcal{T}_{n-1} \subset \mathcal{T}_{n}=\mathscr{T}
  \end{equation}
  by left admissible subcategories~$\mathcal{T}_i$ such that in~$\mathcal{T}_i$, one has~${}^{\perp}\mathcal{T}_{i-1}=\mathcal{N}_i \cong \mathcal{T}_i /\mathcal{T}_{i-1}$.
\end{definition}

In the examples of section \ref{section:applications-constructions} we will mostly consider semi-orthogonal decompositions of a special type, for which the quotients are as simple as possible, namely~$\mathcal{N}_i \cong \derived^\bounded(\catmod/k)$.

\begin{definition}
  An object~$E$ in~$\mathcal{T}$ is \emph{exceptional} if
  \begin{equation}
    \Hom_{\mathcal{T}}(E,E[m])\cong
    \begin{cases}
      k & m=0 \\
      0 & m\neq 0.
    \end{cases}
  \end{equation}
  A sequence of exceptional objects~$(E_1,\ldots,E_l)$ is an \emph{exceptional collection} if
  \begin{equation}
    \Hom_{\mathcal{T}}(E_j,E_i[m])=0 \text{ for } j>i \text{ and any } m
  \end{equation}
  An exceptional collection~$(E_1,\ldots,E_l)$ is called \emph{strong} if in addition
  \begin{equation}
    \Hom_{\mathcal{T}}(E_i,E_j[m])=0 \text{ for all } i,j \text{ and } m \neq 0.
  \end{equation}
  An exceptional collection~$(E_1,\ldots,E_n)$ is called \emph{full} if it generates~$\mathcal{T}$.
\end{definition}

If~$\mathcal{T}$ has a full exceptional collection, then it has a semi-orthogonal decomposition with~$\mathcal{N}_i=\langle E_i \rangle \cong \derived^\bounded(\catmod/k)$ and~$\mathcal{T}_i=\langle E_1,\ldots,E_i\rangle$. We will denote this as
\begin{equation}
  \mathcal{T}=\langle E_1, \ldots, E_n \rangle.
\end{equation}

Recall that a triangulated category is called saturated if all cohomological functors of finite type are representable. Any fully faithful embedding of a saturated triangulated category into a Hom-finite triangulated category is admissible.

From now on let~$X$ denote a smooth projective variety over~$k$. First recall the following result~\cite{baer-tilting-sheaves-representation-theory-algebras,bondal-representation-of-associative-algebras}.

\begin{theorem}
  If~$\derived^\bounded(\coh/X)$ admits a full and strong exceptional collection~$(E_1, \ldots, E_n)$, then the functor
  \begin{equation}
    \RRRHom_X(\bigoplus_{i=1}^n E_i,-)\colon\derived^\bounded(\coh/X) \to \derived^\bounded(\catmod/\End_X(\bigoplus_{i=1}^n E_i))
  \end{equation}
  defines an equivalence of triangulated categories.
\end{theorem}

For~$X$ of dimension~$1$, there is the following theorem of Okawa~\cite{okawa-sod-curve}.

\begin{theorem}
  \label{theorem:okawa}
  The only smooth projective curve~$C$ such that~$\derived^\bounded(\coh/C)$ admits a semi-orthogonal decomposition is~$\mathbb{P}^1$.
\end{theorem}

For dimension~$2$, the results valid for the largest classes of surfaces seem to be due to Hille--Perling and Vial~\cite{hille-perling-tilting-bundles-rational-surfaces-quasi-hereditary-algebras,vial-exceptional-collections-and-the-ns-lattice}.

\begin{theorem}
  A smooth projective rational surface admits a full exceptional collection of line bundles.
\end{theorem}

For higher dimensional varieties, the results are not so general, and we refer to~\cite{kuznetsov-sods-in-algebraic-geometry} for an overview. To construct the exceptional collections that use in section~\ref{section:applications-constructions} we will use the following theorem.

\begin{proposition}[Orlov's blow-up formula]
  \label{proposition:blowup-formula}
  Let~$p$ be a point on~$X$ and consider the blow-up~$\pi\colon\Bl_pX\to X$ with exceptional divisor~$E$. Then
  \begin{equation}
    \derived^\bounded(\coh/\Bl_pX)\cong\langle\pi^*(\derived^\bounded(\coh/X)),\mathcal{O}_E\rangle.
  \end{equation}
\end{proposition}

The motivation for this article is provided by the following result of Orlov~\cite{orlov-nc-schemes}.

\begin{theorem}
 \label{theorem:embedding}
  For a finite dimensional $k$-algebra $A$ of finite global dimension there is an admissible embedding
  \begin{equation}
   \derived^\bounded(\catmod/A) \hookrightarrow \derived^\bounded(\coh/X)
   \end{equation}
  into the derived category of a smooth projective variety $X$.
 \end{theorem}

In fact,~$X$ can be constructed in such a way that it admits a full exceptional collection, so in particular every $A$ as above can be embedded into a triangulated category with a full exceptional collection. Since not every finite-dimensional algebra of finite global dimension admits a full exceptional collection~\cite{happel-fibonacci-algebras}, the existence of such an embedding is non-trivial.

\begin{remark}
  This particular corollary of Orlov's result is already implicitly contained in Iyama's paper~\cite{iyama-finiteness-of-representation-dimension}, which shows that there always exists a quasi-hereditary algebra~$\Lambda$ and an idempotent~$e$ such that~$A=e \Lambda e$, so there is an embedding $\derived^\bounded(\catmod/A) \hookrightarrow \derived^\bounded(\catmod/\Lambda)$. Such an algebra $\Lambda$ always admits a full exceptional collection, provided by the standard modules $\Delta(\lambda)$ coming from the quasi-hereditary structure.
\end{remark}

\section{Embedding finite-dimensional algebras into derived categories of surfaces}
\label{section:obstructions}
We now want to discuss properties of~$A$ that rule out the existence of an embedding
\begin{equation}
  \label{equation:embedding-in-surface}
  \derived^\bounded(\catmod/A) \hookrightarrow \derived^\bounded(\coh/S),
\end{equation}
for some smooth, projective surface~$S$. In section~\ref{section:additive} we consider the (weak) results that one obtains by using noncommutative motives, while we strengthen the results significantly in section~\ref{section:euler-forms} by using extra data.

\subsection{Additive invariants}
\label{section:additive}

Looking for restrictions it is natural to start by checking ``linear'' or additive invariants. An invariant~$I(-)$ of a triangulated category is additive if it is additive with respect to semi-orthogonal decompositions, i.e.
\begin{equation}
  \mathcal{T}=\langle \mathcal{A}, \mathcal{B} \rangle \Rightarrow I(\mathcal{T}) = I(\mathcal{A}) \oplus I(\mathcal{B}).
\end{equation}
Examples of such invariants include algebraic K-theory, nonconnective algebraic K-theory, Hochschild homology, cyclic homology, periodic cyclic homology, negative cyclic homology, topological Hochschild homology and topological cyclic homology, see \cite{tabuada-vandenbergh-noncommutative-motives-of-azumaya-algebras} for some more background.

\begin{lemma}
  Assuming the existence of an embedding~\eqref{equation:embedding-in-surface}, for all of the above invariants~$I(-)$ one has that $I(\derived^\bounded(\mod/A))$ is a direct summand of $I(\derived^\bounded(\coh/S))$.

  \begin{proof}
    Since~$\derived^\bounded(\mod/A)$ has a strong generator, it is saturated~\cite{bondal-vandenbergh}, so~\eqref{equation:embedding-in-surface} is admissible and the claim follows from additivity of $I(-)$.
  \end{proof}
\end{lemma}

However, none of these invariants give particularly interesting information with regards to our question, due to the following (corollary of a) result of Keller~\cite[\S2.5]{keller-invariance-localization-HC-dg-algebras} and Tabuada--Van den Bergh~\cite[corollary~3.20]{tabuada-vandenbergh-noncommutative-motives-of-azumaya-algebras}.

\begin{theorem}
  The additive invariants of a finite-dimensional~$k$\dash algebra of finite global dimension only depend on the number of simple modules.
\end{theorem}


In other words, additive invariants cannot distinguish between algebras with the same number of simple modules, so they are of limited use for our question. Assuming that a specific variety has an additive invariant that is actually computable, this can give a little information though.

\begin{example}
  Using the Hochschild--Kostant--Rosenberg isomorphism for Hochschild homology, it is easy to see that an embedding~\eqref{equation:embedding-in-surface} gives rise to the inequality
  \begin{equation}
    \vert Q_0 \vert=\dim_k\HHHH_0(A)\leq\dim_k\HHHH_0(S)=2+\mathrm{h}^{1,1},
  \end{equation}
  where~$\vert Q_0 \vert$ denotes the number of vertices in the quiver of~$A$, and~$\mathrm{h}^{1,1}$ is the relevant Hodge number.
\end{example}

\begin{remark}
  Of course, all of the above does not depend on~$S$ being a surface and there are obvious generalisations to higher-dimensional varieties.
\end{remark}

\subsection{Quadratic invariants}
\label{section:euler-forms}
We now give two restrictions on the embeddings of derived categories of finite-dimensional algebras in the derived category of a smooth projective surface~$S$. Both results concern the Euler form on the Grothendieck group of such a surface, and are valid for arbitrary surfaces.

Let~$X$ again denote a smooth projective variety of dimension~$n$. Recall that the bilinear Euler form is defined as
\begin{equation}
  \chi\colon\Kzero(X) \times \Kzero(X) \to \mathbb{Z}:(A,B) \mapsto \sum_i (-1)^i \dim_k \Ext^i_X(A,B)
\end{equation}
Moreover one has the natural topological filtration~$\mathrm{F}^\bullet$ on~$\derived^\bounded(\coh/X)$ \cite[Expos\'e~X]{sga6}, where $\mathrm{F}^i\derived^\bounded(\coh/X)$ consists of the complexes of coherent sheaves on~$X$ whose cohomology sheaves have support of codimension at least~$i$.

Denote by~$\serre$ the Serre functor on~$\derived^\bounded(\coh/X)$. Recall that this functor is defined as
\begin{equation}
  \serre(\mathcal{E}^\bullet)=\mathcal{E}^\bullet\otimes\omega_X[n],
\end{equation}
where~$\omega_X$ denotes the canonical line bundle on~$X$. The Serre functor induces an automorphism on~$\Kzero(X)$, which we'll denote by the same symbol, and moreover one has
\begin{equation}
  \label{equation:serre-euler}
  \chi(X,Y)=\chi(Y,\mathbb{S}X).
\end{equation}

By definition,~$\mathrm{F}^i\Kzero(X)$ is the image of the morphism induced by the inclusion~$\mathrm{F}^i\derived^\bounded(\coh/X)\hookrightarrow\derived^\bounded(\coh/X)$. Using this filtration one proves the following result of Suslin \cite[lemma 3.1]{bondal-polishchuk-helices}.
\begin{lemma}
  The operator~$(-1)^n\serre$ is unipotent on~$\Kzero(X)$.
\end{lemma}
To study~$\chi$ using linear algebra, we pass to the numerical Grothendieck group, which is better behaved than the usual Grothendieck group in some respects \cite{mumford-rational-equivalence-on-surfaces}.

\begin{definition}
  The \emph{numerical Grothendieck group}~$\Kzeronum(X)$ is the quotient of~$\Kzero(X)$ by left radical of the Euler form (which by Serre duality agrees with the right radical), i.e.\ one mods out the subgroup of~$\Kzero(X)$ defined by
  \begin{equation}
    \chi(-,\Kzero(X))=\chi(\Kzero(X),-)=0.
  \end{equation}
\end{definition}

For smooth projective varieties this group is always free of finite rank by the Grothendieck--Riemann--Roch theorem, so from now on we will restrict~$\chi$ to~$\Kzeronum(X)$. A closer inspection of the \emph{anti-symmetri\-sation} of the Euler form
\begin{equation}
  \chi^-(A,B)\coloneqq\chi(A,B)-\chi(B,A)=\chi(A,(1-\serre)B) ,
\end{equation}
leads to the first observation. The following result was also proved by Louis de Thanhoffer de V\"olcsey with a different method \cite{louis-phd}.

\begin{theorem}
  \label{theorem:rank}
  For a smooth projective surface~$S$ one has $\rk \chi^{-} \leq 2$.

  \begin{proof}
    By d\'evissage, we can generate~$\Kzero(S)$ by~$[\mathcal{O}_S]$ and classes~$[k_s]$ of skyscrapers for~$s\in S$, and structure sheaves of curves~$[\mathcal{O}_C]$ for all curves~$C$ on~$S$, see~\cite[proposition~0.2.6]{sga6}.

    By Grothendieck--Riemann--Roch the Chern character gives an isomorphism
    \begin{equation}
      \Kzeronum(S)\otimes{\mathbb{Q}}\cong\mathrm{CH}^{\bullet,\mathrm{num}}(S)\otimes\mathbb{Q}
    \end{equation}
    between the numerical Grothendieck group and algebraic cycles modulo numerical equivalence \cite[appendix]{MR2590842}. Hence it suffices to consider a single class~$[k_s]$ as all points are numerically equivalent \cite[\S19.3.5]{fulton-intersection-theory}, and~$\CH^{1,\num}(S)$ is a finitely generated free abelian group of rank~$\rho$ \cite[example~19.3.1]{fulton-intersection-theory}. Observe that under the isomorphism obtained via the Chern character map,~$[k_s]$ is pure of degree~2, $[\mathcal{O}_S]$ is pure of degree~$0$ and~$[\mathcal{O}_C]$ sits in degrees~1 and~2, and the rank of~$\Kzeronum(S)$ is~$\rho+2$.

    We wish to compute the rank of the matrix~$\chi^-$ using this choice of basis. For this we need to know the values of~$\chi^-(\alpha,\beta)$ for~$\alpha\in\CH^i(S)$ and~$\beta\in\CH^j(S)$, with~$i,j\in\{0,1,2\}$. We immediately get that
    \begin{equation}
      \chi^-([k_s],[k_s])=0
    \end{equation}
    and using the presentation
    \begin{equation}
      0\to\mathcal{O}_S(-C)\to\mathcal{O}_S\to\mathcal{O}_C\to0
    \end{equation}
    we get that
    \begin{equation}
      \label{equation:antisym-zero}
      \begin{aligned}
        &\chi^-([\mathcal{O}_C],[k_s]) \\
        &\quad=\chi([\mathcal{O}_S],[k_s])-\chi([k_s],[\mathcal{O}_S])-\chi([\mathcal{O}_S(-C)],[k_s])+\chi([k_s],[\mathcal{O}_S(-C)]) \\
        &\quad=\chi([\mathcal{O}_S],[k_s])-\chi([k_s],[\mathcal{O}_S])-\chi([\mathcal{O}_S],[k_s])+\chi([k_s],[\mathcal{O}_S]). \\
        &\quad=0
      \end{aligned}
    \end{equation}

    For~$i=j=1$ we have that~$C\cdot D=-\chi(\mathcal{O}_C,\mathcal{O}_D)$, hence~$\chi$ is symmetric on this part by the commutativity of the intersection product, therefore it vanishes in the antisymmetric Euler form.

    The (skew-symmetric) matrix one obtains is of the form
    \begin{equation}
      \left(
      \begin{array}{c|c|c}
        0 & \underline{0} & \clubsuit \\ \hline
        \underline{0} & 0\cdot \mathrm{id}_\rho & \underline{\spadesuit} \\ \hline
        \clubsuit & \underline{\spadesuit} & 0
      \end{array}
      \right)
    \end{equation}
    where we order our generators as~$[k_s]$, $[\mathcal{O}_C]$, $[\mathcal{O}_S]$, so there is a block decomposition of a~$(\rho+2)\times(\rho+2)$\dash square matrix with some unknown values, but it is of rank~$\leq 2$ regardless of the unknowns.
  \end{proof}
\end{theorem}

One can also consider the \emph{symmetrised} Euler form
\begin{equation}
  \chi^+(A,B)=\chi(A,B)+\chi(B,A).
\end{equation}
This defines a quadratic form on~$\Kzeronum(X)$, and we can consider its signature, i.e.\ the tuple~$(n_0,n_+,n_-)$ describing the degenerate, positive definite and negative definite part of the form. The forms that we consider are non-degenerate over~$\mathbb{Q}$ by our restriction to~$\Kzeronum(X)$, hence it suffices to specify~$(n_+,n_-)$.

%
%

\begin{theorem}
  \label{theorem:signature}
  Let~$S$ be a smooth projective surface. Then the signature of~$\chi^+$ is~$(\rho,2)$.

  \begin{proof}
    Similar to the proof of theorem~\ref{theorem:rank}, we consider the (ordered) basis $[k_s]$, $[\mathcal{O}_S]$, $[\mathcal{O}_C]_C$. The Hodge index theorem \cite[example~19.3.1]{fulton-intersection-theory} says that the signature for the intersection product on the curves is~$(1,\rho-1)$. Via the equality~$C\cdot D=-\chi(\mathcal{O}_C,\mathcal{O}_D)$, the subspace spanned by the $[\mathcal{O}_C]$ has signature~$(\rho-1,1)$ for $\chi^+$.

    Because~$\chi^+([k_s],[k_s])=0$, $\chi^+([\mathcal{O}_S],[k_s])=2$ using Serre duality and~$\chi^+([\mathcal{O}_S],[\mathcal{O}_S])$ is twice the Euler characteristic of~$S$, the subspace spanned by $[k_s]$ and $[\mathcal{O}_S]$ is a hyperbolic plane, so it has signature~$(1,1)$.

    For the other terms we get that~$\chi^+([\mathcal{O}_C],[k_s])=0$ just as in \eqref{equation:antisym-zero}, whilst~$\chi^+([\mathcal{O}_C],[\mathcal{O}_S])$ can be arbitrary. Therefore, after a base change, the matrix of~$\chi^+$ has the form
    \begin{equation}
      \label{equation:signature-decomposition}
      \left(
      \begin{array}{c|c|c}
        0 & 2 & \underline{0} \\ \hline
        2 & 2\chi(S) & \underline{*}^{\mathrm{t}} \\ \hline
        \underline{0} & \underline{*} & D
      \end{array}
    \right)
  \end{equation}
  where~$D=\text{diag}(-1,1,\ldots,1)$ is a diagonal matrix of size~$\rho$. In particular, the decomposition is not orthogonal.

  If~$\rho=1$ then the sign of the determinant of the matrix is positive, but there is at least one negative eigenvalue coming from~$D$. So the signature must be~$(1,2)$.

  If~$\rho\geq 2$, denote by~$W$ the subspace spanned by the last~$\rho-2$ basis vectors for the choice of basis as in \eqref{equation:signature-decomposition}. The subspace spanned by~$(1,0,\ldots,0)$ and~$(0,0,1,1,0,\ldots,0)$ is a totally isotropic subspace in~$W^{\perp}$. Since~$\chi^+$ is non-degenerate~$W^{\perp}$ contains an orthogonal sum of two hyperbolic planes. Now~$\dim \Kzeronum(S) \otimes \mathbb{Q}= \dim W + \dim W^{\perp}$, so~$W^{\perp}\sim2\mathbb{H}$. We conclude that the signature of~$\chi^+$ can be computed as
  \begin{equation}
    (\rho-2,0) + 2(1,1)=(\rho,2).
    \label{}
  \end{equation}
  \end{proof}
\end{theorem}

We now apply these observations to give restrictions on embeddings of (bounded derived categories of) finite-dimensional algebras into (bounded derived categories of) smooth projective surfaces. Let~$A$ denote a basic finite-dimensional~$k$\dash algebra of finite global dimension with~$n$ simple modules. In that case, there is a well-defined Euler form given by
\begin{equation}
  \chi\colon\Kzero(A) \times \Kzero(A) \to \mathbb{Z}:(X,Y) \mapsto \sum_i (-1)^i \dim_k \Ext^i_A(X,Y).
\end{equation}
Since the indecomposable projective modules and the simple modules form dual bases, this bilinear form is non-degenerate. Also submatrices cannot increase in rank and signatures behave well under restriction, so the following corollaries are clear.
\begin{corollary}
  \label{corollary:rank}
  Given a smooth projective surface~$S$ and an embedding
  \begin{equation}
    \derived^\bounded(\catmod/A) \hookrightarrow \derived^\bounded(\coh/X),
  \end{equation}
  the rank of~$\chi^-_A$ is~$\leq 2$.
\end{corollary}
\begin{corollary}
  \label{corollary:signature}
  Given a smooth projective surface~$S$ and an embedding
  \begin{equation}
    \derived^\bounded(\catmod/A) \hookrightarrow \derived^\bounded(\coh/X),
  \end{equation}
  then $\chi_{A}^+$ does not admit a~$3$\dash dimensional negative definite subspace.
\end{corollary}

\section{Embedding hereditary algebras}
\label{section:applications-constructions}
In this section we show how the established criteria can be applied to restrict embeddings as in theorem~\ref{theorem:embedding} for hereditary algebras, as these have a particularly nice description for their derived categories. Observe that by using tilting theory it is possible to find finite-dimensional algebras of global dimension~$\geq 2$ which are derived equivalent to path algebras. Hence all results in this section are also valid for iterated tilted algebras.

The extreme case of the embedding being an equivalence has a particularly easy answer. It is elementary that~$A$ being semisimple implies~$X$ is a union of points. In global dimension~$1$ there is the following easy folklore result.

\begin{proposition}
  \label{proposition:hereditary-endomorphism-algebra-tilting-module}
  If~$A=kQ$ is hereditary (and not semisimple), and
  \begin{equation}
    \derived^\bounded(\coh/X) \cong \derived^\bounded(\catmod/A)
  \end{equation}
  is an equivalence, then~$X\cong\mathbb{P}^1$ and~$A \cong k\kronecker_2$, the path algebra of the Kronecker quiver.

  \begin{proof}
    The description of~$\mathbb{P}^1$ is standard~\cite{beilinson-coherent-linear-algebra}. To see that this is the only variety with this property consider the skyscraper sheaves~$k_x$, which are indecomposable objects (or more precisely, they are point objects).

    A triangle equivalence sends these to indecomposable objects of~$\derived^\bounded(\catmod/kQ)$, which correspond to the indecomposable modules up to a shift since every object therein is formal. Now by Serre duality
    \begin{equation}
      \begin{aligned}
        \Ext_X^d(k_x,k_x)&\cong \Hom_X(k_x[d],k_x \otimes \omega_X[\dim X])^\vee \\
        &\cong \Hom_X(k_x,k_x[\dim X -d])^\vee,
      \end{aligned}
    \end{equation}
    and since a hereditary algebra is of global dimension~$1$,~$X$ has to be a curve.

    By by Serre duality we have that for every object~$E \in \derived^\bounded(\coh/X)$
    \begin{equation}
      \Hom_X(E,E[1]) \cong \Hom_X(E,E\otimes \omega_X)^\vee \neq 0.
    \end{equation}
    In particular,~$\derived^\bounded(\coh/X)$ only contains exceptional objects if~$X \cong \mathbb{P}^1$, otherwise a nonzero section of~$\omega_X$ gives a nonzero morphism.

    Any hereditary algebra derived equivalent to~$X$ is thus derived equivalent to~$k\kronecker_2$. It is known, see~\cite[\S 4.8]{happel-on-the-derived-category-of-a-finite-dimensional-algebra}, that any derived equivalence between basic hereditary algebras is given by a sequence of sink or source reflections, so there is no possibility other than~$\kronecker_2$.
  \end{proof}
\end{proposition}

Having settled the case of an equivalence assume there is an embedding
 \begin{equation}
   \label{equation:embedding-path-algebra}
   \derived^\bounded(\catmod/kQ) \hookrightarrow \derived^\bounded(\coh/S),
 \end{equation}
for an acyclic quiver~$Q$ and a smooth projective surface~$S$, i.e.\ we are in the situation of a strong (but not full) exceptional collection without relations in the composition law. The following questions regarding the structure of possible quivers~$Q$ come to mind:
\begin{enumerate}[label=(Q\arabic*)]
  \item\label{enumerate:quiver-question-1} Is there a bound on the number of vertices of~$Q$?
  \item\label{enumerate:quiver-question-2} Is there a bound on the number of arrows of~$Q$?
  \item\label{enumerate:quiver-question-3} Is there a bound on the number of paths in~$Q$?
  \item\label{enumerate:quiver-question-4} Is there a bound on the length of paths in~$Q$?
  \item\label{enumerate:quiver-question-5} Is it possible to embed any quiver on~$3$ vertices?
\end{enumerate}

The remainder of this section is dedicated to answering these questions. In the following proposition we provide some explicit embeddings for well known families of quivers.

\begin{proposition}
  Let~$A=kQ$ be the path algebra of an acyclic quiver~$Q$.

  \begin{enumerate}
    \item If~$A$ is of finite type or tame, i.e.~$Q$ is a Dynkin or Euclidean quiver, then an embedding~\eqref{equation:embedding-path-algebra} exists if and only if
      \begin{equation*}
        Q=\An_1,\An_2,\An_3,\Dn_4,\tilde{\An}_1 \text{ or } \tilde{\An}_2.
      \end{equation*}

    \item If~$Q$ occurs in the following families of quivers:
      \begin{equation}
        \begin{gathered}
          \kronecker_n\colon
          \begin{tikzpicture}[scale = 0.7, baseline = -.5ex]
            \node (1)                {};
            \node (2) [right of = 1] {};

            \draw (1) circle (2pt);
            \draw (2) circle (2pt);

            \draw[->] (1) to [bend left = 45] node [above] {$1$} (2);
            \draw[->] (1) to [bend right = 45] node [below] {$n$} (2);
            \draw[color = white] (1) to node [color = black, rotate = 90] {$\dots$} (2);
          \end{tikzpicture}, \quad\quad
          \Sn_n\colon
          \begin{tikzpicture}[scale = 0.7, baseline = 0pt]
            \node (0)                {};
            \node (2) [right of = 0] {\vdots};
            \node (1) [above of = 2] {};
            \node (n) [below of = 2] {};

            \draw (0) circle (2pt) node[above] {$0$};
            \draw (1) circle (2pt) node[above] {$1$};
            \draw (n) circle (2pt) node[below] {$n$};

            \draw[->] (0) -- (1);
            \draw[->] (0) -- (n);
          \end{tikzpicture}
        \end{gathered}
      \end{equation}
      then an embedding~\eqref{equation:embedding-path-algebra} exists.
  \end{enumerate}

  \begin{proof}
    For part~(1), it suffices to compute the matrices for the anti-symmetric Euler forms which can be obtained from the Cartan matrices, and observe that starting from~$\An_4$, $\Dn_5$ and~$\tilde{\An}_3$ the rank is bounded below by~4 as there will be a submatrix of rank~4\ in each of these coming from the smallest cases~$\An_4$, $\Dn_5$ and~$\tilde{\An}_3$. The exceptional types~$\En_{6,7,8}$ or~$\tilde{\En}_{6,7,8}$ have corresponding ranks~$6,6,8,6,6$ and~$8$.

    For the 5~cases that are not ruled out by this restriction on the anti-symmetric Euler form, and the infinite families in part~(2), there are explicit embeddings.

    Let~$n=2m$, then one can embed~$\kronecker_{2m}$ by considering~$\mathcal{O}_{\mathbb{P}^1\times\mathbb{P}^1}$ and~$\mathcal{O}_{\mathbb{P}^1\times\mathbb{P}^1}(1,m-1)$ on~$\mathbb{P}^1 \times \mathbb{P}^1$.

    Let~$n=2m-1$, then one can embed~$\kronecker_{2m-1}$ by considering~$\mathcal{O}_{\Bl_1\mathbb{P}^2}$ and $\mathcal{O}_{\Bl_1\mathbb{P}^2}(E+mF)$, on the blow-up of~$\mathbb{P}^2$ in a point~$p$. Here, as usual,~$E$ denotes the divisor associated to the~$-1$\dash curve and~$F$ the one associated to the strict transform of any line in~$\mathbb{P}^2$ through~$p$.

    For the family~$\Sn_n$, by Orlov's blow-up formula in~proposition~\ref{proposition:blowup-formula} we obtain a semi-orthogonal decomposition
    \begin{equation}
      \derived^\bounded(\coh/\Bl_1\mathbb{P}^2)=\langle\pi^*(\derived^\bounded(\coh/\mathbb{P}^2)),\mathcal{O}_E\rangle
    \end{equation}
    where~$\pi\colon\Bl_1\mathbb{P}^2\to\mathbb{P}^2$ is the blow-up morphism, and~$\mathcal{O}_E$ is the structure sheaf of the exceptional divisor. The blow-up locus is denoted~$p$. Consider the exceptional line bundle~$\mathcal{O}_{\mathbb{P}^2}$ on~$\mathbb{P}^2$, then one checks by adjunction
    \begin{equation}
      \Hom_{\derived^\bounded(\coh/\Bl_1\mathbb{P}^2)}(\pi^*(\mathcal{O}_{\mathbb{P}^2}),\mathcal{O}_E)\cong\Hom_{\derived^\bounded(\coh/\mathbb{P}^2)}(\mathcal{O}_{\mathbb{P}^2},k_p)
    \end{equation}
    that the exceptional pair~$(\pi^*(\mathcal{O}_{\mathbb{P}^2}),\mathcal{O}_E$ has endomorphism ring~$k\Sn_1$. Using the blow-up formula inductively, this gives a realisation of~$k\Sn_n$ using~$\Bl_n\mathbb{P}^2$.

    By the identifications~$\An_2=\Sn_1$, $\An_3=\Sn_2$ (using reflection), $\An_4=\Sn_3$ (using reflection) and~$\kronecker_2=\tilde{\An}_1$ the only remaining quiver is~$\tilde{\An}_2$, and for this one we use some elementary toric geometry \cite{MR2810322}. The variety ~$\Bl_2\mathbb{P}^2$ can be represented by the fan
    \begin{equation}
      \label{equation:fan}
      \begin{tikzpicture}[scale = .7, baseline = 0em]
        \draw[->] (0, 0) -- (1, 0) node [right] {1};
        \draw[->] (0, 0) -- (0, 1) node [above] {2};
        \draw[->] (0, 0) -- (-1, 0) node [left] {3};
        \draw[->] (0, 0) -- (-1, -1) node [left] {4};
        \draw[->] (0, 0) -- (0, -1) node [below] {5};
      \end{tikzpicture}
    \end{equation}
    As basis for~$\Pic(\Bl_2\mathbb{P}^2)$, we choose the first three torus-invariant divisors~$D_1,D_2$ and~$D_3$. It can be computed that
    \begin{equation}
      (\mathcal{O}_{\Bl_2\mathbb{P}^2},\mathcal{O}_{\Bl_2\mathbb{P}^2}(D_1-D_3),\mathcal{O}_{\Bl_2\mathbb{P}^2}(D_1))
    \end{equation}
    has the desired structure.
  \end{proof}
\end{proposition}

\begin{remark}
  Of course, there are many alternatives to the above embeddings. The Kronecker quivers~$\kronecker_n$ for example can also be embedded using~$\mathcal{O}(E)$ and~$\mathcal{O}(E+nF)$ on~$\mathbf{F}_n=\PProj(\mathcal{O}_{\mathbb{P}^1} \oplus \mathcal{O}_{\mathbb{P}^1}(-n))$, the~$n$th Hirzebruch surface.
\end{remark}

Using this proposition, \ref{enumerate:quiver-question-1} and~\ref{enumerate:quiver-question-2} have a negative answer. For~\ref{enumerate:quiver-question-3} the answer is also no, since one can always reflect the~$\Sn_n$\dash quiver in some non-zero vertex.

The questions~\ref{enumerate:quiver-question-4} and~\ref{enumerate:quiver-question-5} are more subtle. Let us say a quiver~$Q'$ is \emph{forbidden} if the rank of~$\chi^-$ is strictly greater than~$2$.

\begin{lemma}
  If a quiver~$Q$ contains~a forbidden quiver~$Q'$ as a full subquiver, then it cannot be embedded into a smooth projective surface.

  \begin{proof}
   The fullness ensures that the~$\chi^-$ matrix of~$Q'$ occurs as a block in that of~$Q$ (for the basis of simples for example), so~$\rk(\chi^-_Q) >2$ and the quiver cannot be embedded.
  \end{proof}
\end{lemma}

Then question \ref{enumerate:quiver-question-4} about path length can be partially answered by plugging~$\An_4$ into this lemma, as we know from proposition~\ref{corollary:rank} that~$\An_4$ cannot be embedded into the derived category of a surface. Observe that~$\An_4$ does satisfy the condition on the negative definite subspaces for~$\chi^+$, as in proposition~\ref{corollary:signature}.

However, if~$\An_4$ occurs as a non-full subquiver of~$Q$, it is not clear what happens. The following example shows that in some cases one \emph{can} have an embedding.

\begin{example}
  Consider the following quiver:
  \begin{equation}
    \mathsf{Q}\colon
    \begin{tikzpicture}[baseline = -.5ex]
      \node (1)                {};
      \node (2) [right of = 1] {};
      \node (3) [right of = 2] {};
      \node (4) [right of = 3] {};

      \draw (1) circle (2pt);
      \draw (2) circle (2pt);
      \draw (3) circle (2pt);
      \draw (4) circle (2pt);

      \draw[->] (1) to [bend left = 45] (3);
      \draw[->] (1) to (2);
      \draw[->] (2) to [bend right = 45] (4);
      \draw[->] (2) to (3);
      \draw[->] (3) to (4);
    \end{tikzpicture}
  \end{equation}
  It has~$\rk(\chi^-)=2$ and contains~$\An_4$ but only as a non-full subquiver, so it does not satisfy the condition for the previous proposition. In fact, it can be embedded into~$\Bl_2 \mathbb{P}^2$ by extending the strong exceptional collection we already had for~$\tilde{\An}_2$. In terms of the fan mentioned in~\eqref{equation:fan}, the collection
  \begin{equation}
    (\mathcal{O}_{\Bl_2 \mathbb{P}^2},\mathcal{O}_{\Bl_2 \mathbb{P}^2}(D_1-D_3), \mathcal{O}_{\Bl_2 \mathbb{P}^2}(D_1),\mathcal{O}_{\Bl_2 \mathbb{P}^2}(D_2))
  \end{equation}
  can be checked to be a strong exceptional collection with the desired endomorphism ring, after reflecting in the vertex corresponding to~$\mathcal{O}_{\Bl_2 \mathbb{P}^2}$.
\end{example}

It is also possible to find an example of a quiver that satisfies proposition~\ref{corollary:rank} but violates proposition~\ref{corollary:signature}.
\begin{example}
  Consider an acyclic quiver on~5 vertices~$v_1,\dotsc,v_n$, whose~$\chi$ is given by the matrix
  \begin{equation}
    \begin{pmatrix}
      1 & 2 & 4 & 3 & 0 \\
      0 & 1 & 4 & 5 & 2\\
      0 & 0 & 1 & 4 & 4 \\
      0 & 0 & 0 & 1 & 3 \\
      0 & 0 & 0 & 0 & 1
    \end{pmatrix}.
  \end{equation}
  It is straightforward to check that~$\rk\chi^-=2$, but~$\chi^+$ has a negative-definite subspace of dimension~3.
\end{example}

The reason for posing question \ref{enumerate:quiver-question-5} is that any skew-symmetric~$3 \times 3$\dash matrix has rank~$\leq 2$, and moreover it cannot have a~$3$\dash dimensional negative definite subspace, since one can always look at the projective indecomposables which yield a nonzero positive definite subspace. Also, since every quiver on~$2$ vertices can be embedded, \ref{enumerate:quiver-question-5} naturally arises as the next case.

Such a quiver can be presented as
\begin{equation}
  \label{equation:3-vertex-quiver}
  Q_{a,b,c}\colon
  \begin{tikzpicture}[baseline = -20pt, node distance = 1.5cm]
    \node (1)                {};
    \node (2) [right of = 1] {};
    \node (3) [below of = 2] {};

    \draw (1) circle (2pt) node[above]  {};
    \draw (2) circle (2pt) node[right]  {};
    \draw (3) circle (2pt) node[above]  {};

    \draw[->, bend left = 25]  (1) edge (2);
    \draw[->, bend right = 25] (1) edge (2);
    \draw[opacity = 0]         (1) edge node[color = black, opacity = 1] {\raisebox{.2cm}{\tiny$\vdots\raisebox{.1cm}{$a$}$}} (2);

    \draw[->, bend left = 25]  (2) edge (3);
    \draw[->, bend right = 25] (2) edge (3);
    \draw[opacity = 0]         (2) edge node[color = black, opacity = 1] {\raisebox{.2cm}{\tiny$\vdots\raisebox{.1cm}{$b$}$}} (3);

    \draw[->, bend left = 20]  (1) edge (3);
    \draw[->, bend right = 20] (1) edge (3);
    \draw[opacity = 0]         (1) edge node[color = black, opacity = 1] {\raisebox{.2cm}{\tiny$\vdots\raisebox{.1cm}{$c$}$}} (3);
  \end{tikzpicture}
\end{equation}

To make computations feasible we will only consider line bundles on rational surfaces, i.e.\ iterated blow-ups of the minimal rational surfaces~$\mathbb{P}^2$ and~$\mathbf{F}_n$ for~$n\neq 1$. In this context we can apply techniques based on Riemann--Roch arithmetic which are for instance also used in \cite{hille-perling-exceptional-sequences-of-invertible-sheaves}.

\begin{theorem}
  The values~$(a,b,c)$ as in \eqref{equation:3-vertex-quiver} for which there exists a rational surface~$S$ and an embedding
  \begin{equation}
	  \label{equation:embedding-3-vertex}
	  \derived^\bounded(\catmod/kQ_{a,b,c}) \hookrightarrow \derived^\bounded(\coh/S)
	\end{equation}
  given by a strong exceptional collection of line bundles are
  \begin{equation}
    \begin{gathered}
      \left\{ (0,n,n)\mid n\in\mathbb{N} \right\}\cup\left\{ (n,0,n)\mid n\in\mathbb{N} \right\}, \\
      \left\{ (1,n,1)\mid n\in\mathbb{N} \right\}\cup\left\{ (n,1,1)\mid n\in\mathbb{N} \right\}, \\
      \left\{ (2,2,0) \right\}.
    \end{gathered}
  \end{equation}

  \begin{proof}
    By twisting we have that the exceptional collection is of the form
    \begin{equation}
      \langle \mathcal{O}_S, \mathcal{O}_S(D),\mathcal{O}_S(E) \rangle,
    \end{equation}
    for divisors~$D$ and~$E$.

    We first claim that if an embedding~\eqref{equation:embedding-3-vertex} exists, then
    \begin{equation}
      \label{equation:3-vertex-constraints}
      a+b=ab+c.
    \end{equation}
    To see this, first note that on a rational surface~$\chi(\mathcal{O}_S)=1$. From exceptionality we get~$0=\chi(\mathcal{O}_S(D),\mathcal{O})=\chi(-D)$, and similarly~$\chi(-E)=\chi(D-E)=0$. By Riemann--Roch we obtain
	  \begin{equation}
	    \chi(D)=\chi(\mathcal{O}_S) + \frac{1}{2}(D^2 -K_S \cdot D),
	  \end{equation}
    for any divisor~$D$. By anti-symmetrising this equation for our divisors~$D, E$ and~$E-D$ respectively and plugging in the zeroes we found above, we get
	  \begin{equation}
	    \chi(D)+\chi(E-D)=\chi(E).
	  \end{equation}
    By strong exceptionality, $\mathrm{h}^i(\mathcal{O}_S(D))=\mathrm{h}^i(\mathcal{O}_S(E))=\mathrm{h}^i(\mathcal{O}_S(E-D))=0$ for all~$i>0$, so~$\chi(D)=\mathrm{h}^0(\mathcal{O}_S(D))$ and similarly for~$E$ and~$E-D$, so we finally find
	  \begin{equation}
      \mathrm{h}^0(\mathcal{O}_S(D))+\mathrm{h}^0(\mathcal{O}_S(E-D))=\mathrm{h}^0(\mathcal{O}_S(E)) \Rightarrow a+b=ab+c.
	  \end{equation}

    Solving this equation yields the solutions listed in the statement. We now give a construction for each of these cases. All of these can be realised on the del Pezzo surface~$\Bl_3\mathbb{P}^2$ of degree~6, which will be represented by the fan
    \begin{equation}
      \begin{tikzpicture}[scale = .7, baseline = 0]
        \draw[->] (0,0) -- (1,0) node [right] {1};
        \draw[->] (0,0) -- (1,1) node [above] {2};
        \draw[->] (0,0) -- (0,1) node [above] {3};
        \draw[->] (0,0) -- (-1,0) node [left] {4};
        \draw[->] (0,0) -- (-1,-1) node [below] {5};
        \draw[->] (0,0) -- (0,-1) node [below] {6};
      \end{tikzpicture}
    \end{equation}
    and as a basis for~$\Pic(X)$ we choose the torus-invariant divisors~$D_1,\dotsc,D_4$.

    \begin{table}
      \centering
      \footnotesize
      \begin{tabular}{clll}
        \toprule
        $(a,b,c)$    &          & $D$                 & $E$ \\\midrule
        $(0,n,n)$    & $n=2m$   & $D_2-D_4$           & $(m-1)D_1+mD_2+D_3$ \\
                     & $n=2m+1$ & $D_1-D_4$           & $D_1+D_2+mD_3+(m-1)D_4$ \\\addlinespace
        $(n,0,n)$    & $n=2m$   & $D_2+mD_3+(m-1)D_4$ & $D_1+D_2+(m-1)D_3+(m-1)D_4$ \\
                     & $n=2m+1$ & $D_2+mD_3+mD_4$     & $D_1+D_2+mD_3+(m-1)D_4$ \\\addlinespace
        $(1,n,1)$    & $n=2m$   & $D_1+D_2-D_4$       & $mD_1+mD_2+D_3$ \\
                     & $n=2m+1$ & $D_4$               & $mD_1+mD_2+D_3+D_4$ \\\addlinespace
        $(n,1,1)$    & $n=2m$   & $(m-1)D_1+mD_2+D_3$ & $(m-1)D_1+mD_2+D_3+D_4$ \\
                     & $n=2m+1$ & $D_2+mD_3+mD_4$     & $D_1+D_2+mD_3+mD_4$ \\\bottomrule
      \end{tabular}
      \caption{Divisors for embeddings of 3-vertex quivers in~$\derived^\bounded(\coh/\Bl_3\mathbb{P}^2)$}
      \label{table:divisors}
    \end{table}

    For the families of solutions we give a possible choice of divisors in table~\ref{table:divisors}, and code to check these results in the appendix.

    For the isolated case~$(2,2,0)$ it suffices to take~$D=L_1$ and~$E=L_1+L_2$ on~$\mathbb{P}^1\times\mathbb{P}^1$, where the~$L_i$ form the basis of~$\Pic(\mathbb{P}^1\times\mathbb{P}^1)$ given by the two rulings.
  \end{proof}
\end{theorem}

\begin{remark}
  On a rational surface any exceptional sequence can be mutated into one consisting of rank one objects, and the previous result also holds in this greater generality, since we only used numerical computations to obtain~\eqref{equation:3-vertex-constraints}. We leave open the question whether considering more general exceptional vector bundles (or coherent sheaves) on more general surfaces yields more examples.
\end{remark}

\appendix
\section{Code for table \ref{table:divisors}}
The following Sage-code checks the results from table~\ref{table:divisors}.
{
\small
\begin{verbatim}
X = toric_varieties.dP6()
G = X.rational_class_group()

def checkQuiver((a, b, c), D, E):
  print "(a,b,c) = (%d,%d,%d) using D = %s and E = %s" % (a, b, c, D, E)

  D, E = G(D).lift(), G(E).lift()

  # checking exceptional pairs
  assert (-D).cohomology(dim=True).is_zero()
  assert (-E).cohomology(dim=True).is_zero()
  assert (D-E).cohomology(dim=True).is_zero()

  # structure of the quiver
  assert tuple(D.cohomology(dim=True)) == (a, 0, 0)
  assert tuple((E-D).cohomology(dim=True)) == (b, 0, 0)
  assert tuple(E.cohomology(dim=True)) == (a*b + c, 0, 0)

for m in range(0, 10):
  checkQuiver((0, 2*m, 2*m),     [0, 1, 0, -1],  [m-1, m, 1, 0])
  checkQuiver((0, 2*m+1, 2*m+1), [1, 0, 0, -1],  [1, 1, m, m-1])
  checkQuiver((2*m, 0, 2*m),     [0, 1, m, m-1], [1, 1, m-1, m-1])
  checkQuiver((2*m+1, 0, 2*m+1), [0, 1, m, m],   [1, 1, m, m-1])
  checkQuiver((1, 2*m, 1),       [1, 1, 0, -1],  [m, m, 1, 0])
  checkQuiver((1, 2*m+1, 1),     [0, 0, 0, 1],   [m, m, 1, 1])
  checkQuiver((2*m, 1, 1),       [m-1, m, 1, 0], [m-1, m, 1, 1])
  checkQuiver((2*m+1, 1, 1),     [0, 1, m, m],   [1, 1, m, m])
\end{verbatim}
}

\bibliography{biblio}

\providecommand{\bysame}{\leavevmode\hbox to3em{\hrulefill}\thinspace}
\providecommand{\MR}{\relax\ifhmode\unskip\space\fi MR }
\providecommand{\MRhref}[2]{%
  \href{http://www.ams.org/mathscinet-getitem?mr=#1}{#2}
}
\providecommand{\href}[2]{#2}
\begin{thebibliography}{10}

\bibitem{baer-tilting-sheaves-representation-theory-algebras}
Dagmar Baer, \emph{Tilting sheaves in representation theory of algebras},
  \textbf{60} (1988), 323--347.

\bibitem{beilinson-coherent-linear-algebra}
Alexander Be{\u\i}linson, \emph{Coherent sheaves on $\mathbb{P}^n$ and problems
  of linear algebra},  \textbf{12} (1978), no.~3, 214--216.

\bibitem{sga6}
Pierre Berthelot, Alexander Grothendieck, and Luc Illusie (eds.),
  \emph{S\'emi\-naire de g\'eom\'etrie alg\'ebrique du bois marie --- 1966--67
  --- th\'eorie des intersections et th\'eor\`eme de riemann--roch --- (sga6)},
  Lecture notes in mathematics, no. 225, Springer-Verlag, 1971.

\bibitem{bondal-representation-of-associative-algebras}
Alexei Bondal, \emph{Representation of associative algebras and coherent
  sheaves},  \textbf{34} (1990), no.~1, 23--42.

\bibitem{bondal-polishchuk-helices}
Alexei Bondal and Alexander Polishchuk, \emph{Homological properties of
  associative algebras: the method of helices},  \textbf{42} (1994), no.~2,
  219--260.

\bibitem{bondal-vandenbergh}
Alexei Bondal and Michel Van~den Bergh, \emph{Generators and representability
  of functors in commutative and noncommutative geometry},  \textbf{3} (2003),
  1--36 (english).

\bibitem{MR2810322}
David~A. Cox, John~B. Little, and Henry~K. Schenck, \emph{Toric varieties},
  Graduate Studies in Mathematics, vol. 124, American Mathematical Society,
  Providence, RI, 2011.

\bibitem{louis-phd}
Louis {de Thanhoffer de V\"olcsey}, \emph{Non-commutative projective geometry
  and calabi--yau algebras}, Ph.D. thesis, 2015.

\bibitem{fulton-intersection-theory}
William Fulton, \emph{Intersection theory}, 2 ed., Ergebnisse der Mathematik
  und ihrer Grenzgebiete, 3.\ Folge, Springer, 1998.

\bibitem{gorchinskiy-orlov-geometric-phantom-categories}
Sergey Gorchinskiy and Dmitri Orlov, \emph{Geometric phantom categories},
  Publications Math\'ematiques, Institut Hautes \'Etudes Scientifiques
  \textbf{117} (2013), 329--349.

\bibitem{happel-on-the-derived-category-of-a-finite-dimensional-algebra}
Dieter Happel, \emph{On the derived category of a finite-dimensional algebra},
  \textbf{62} (1987), 339--389.

\bibitem{happel-fibonacci-algebras}
\bysame, \emph{A family of algebras with two simple modules and fibonacci
  numbers},  \textbf{57} (1991), 133--139.

\bibitem{hille-perling-exceptional-sequences-of-invertible-sheaves}
Lutz Hille and Markus Perling, \emph{Exceptional sequences of invertible
  sheaves on rational surfaces},  \textbf{147} (2011), 1230--1280.

\bibitem{hille-perling-tilting-bundles-rational-surfaces-quasi-hereditary-algebras}
\bysame, \emph{Tilting bundles on rational surfaces and quasi-hereditary
  algebras},  (2011).

\bibitem{iyama-finiteness-of-representation-dimension}
Osamu Iyama, \emph{Finiteness of representation dimension},  \textbf{131}
  (2002), no.~4, 1011--1014.

\bibitem{keller-invariance-localization-HC-dg-algebras}
Bernhard Keller, \emph{Invariance and localization for cyclic homology of dg
  algebras},  \textbf{123} (1998), 223--273.

\bibitem{MR2590842}
A.~G. Kuznetsov, \emph{Derived categories of {F}ano threefolds}, Tr. Mat. Inst.
  Steklova \textbf{264} (2009), 116--128.

\bibitem{kuznetsov-sods-in-algebraic-geometry}
Alexander Kuznetsov, \emph{Semiorthogonal decompositions in algebraic
  geometry}, Proceedings of the ICM 2014, 2014.

\bibitem{mumford-rational-equivalence-on-surfaces}
David Mumford, \emph{Rational equivalence of 0-cycles on surfaces},  \textbf{9}
  (1969), 1969.

\bibitem{okawa-sod-curve}
Shinnosuke Okawa, \emph{Semi-orthogonal decomposability of the derived category
  of a curve},  \textbf{228} (2011), 2869--2873.

\bibitem{orlov-nc-schemes}
Dmitri Orlov, \emph{Smooth and proper noncommutative schemes and gluing of dg
  categories},  (2014).

\bibitem{orlov-geometric-realization}
\bysame, \emph{Geometric realizations of quiver algebras}, 2015.

\bibitem{tabuada-vandenbergh-noncommutative-motives-of-azumaya-algebras}
Gon\c{c}alo Tabuada and Michel Van~den Bergh, \emph{Noncommutative motives of
  azumaya algebras},  (2014), 1--25.

\bibitem{vial-exceptional-collections-and-the-ns-lattice}
Charles Vial, \emph{Exceptional collections, and the n\'eron--severi lattice
  for surfaces}, 2015.

\end{thebibliography}
\bibliographystyle{amsplain} 

\end{document}